\documentclass[12pt,a4paper]{article}   

\input amssym
 
\catcode`\@=11

\font\tenmsa=msam10 \font\sevenmsa=msam7 \font\fivemsa=msam5
\font\tenmsb=msbm10 at 12pt\font\sevenmsb=msbm7 at 9pt\font\fivemsb=msbm5 at 7pt
\newfam\msafam \newfam\msbfam
\textfont\msafam=\tenmsa  \scriptfont\msafam=\sevenmsa \scriptscriptfont\msafam=
\fivemsa
\textfont\msbfam=\tenmsb  \scriptfont\msbfam=\sevenmsb \scriptscriptfont\msbfam=
\fivemsb

\font\teneufm=eufm10 \font\seveneufm=eufm7 \font\fiveeufm=eufm5
\newfam\eufmfam 
\textfont\eufmfam=\teneufm \scriptfont\eufmfam=\seveneufm\scriptscriptfont
\eufmfam=\fiveeufm

\def\hexnumber@#1{\ifcase#1 0\or1\or2\or3\or4\or5\or6\or7\or8\or9\or
 A\or B\or C\or D\or E\or F\fi }

\edef\msa@{\hexnumber@\msafam} \edef\msb@{\hexnumber@\msbfam}

\mathchardef\Box="0\msa@03

\catcode`\@=12

\newtheorem{thr}{Theorem}[section]
\newtheorem{rem}[thr]{Remark}
\newtheorem{cor}[thr]{Corollary}
\newtheorem{prop}[thr]{Proposition}
\newtheorem{defi}[thr]{Definition}
\newtheorem{lemma}[thr]{Lemma}

\def\RR{\Bbb R}
\def\NN{\Bbb N}
\def\CC{\Bbb C}
\def\fst{f_2^{*}}
\def\luw{\Lambda^p_2(w)}

\def\Rdm{{\RR}^2_+}

\begin{document}
\title{Multidimensional rearrangement \\ and Lorentz spaces}
\author{ {\bf Sorina Barza} \\ Department of Mathematics, Karlstad University\\ S-65188 
Karlstad, SWEDEN \\E-mail: {\tt sorina.barza@kau.se} \and {\bf  Lars-Erik Persson}\\ Department
of Mathematics, Lule\aa\, University\\
 S-97187  Lule\aa, SWEDEN\\E-mail: {\tt larserik@sm.luth.se}  \and {\bf Javier
Soria\thanks{Partially supported
 by DGICYT PB97-0986 and CIRIT 1999SGR00061} }\\Dept.\
Appl.\ Math.\  and Analysis, University of  Barcelona\\
E-08071 Barcelona, SPAIN\\E-mail: {\tt soria@mat.ub.es}}
\date{} 
\maketitle
\begin{center}{\bf Abstract}
\end{center}
We define a multidimensional rearrangement, which is related to classical inequalities for
functions that are monotone in each variable. We prove the main measure theoretical results of
the new theory and characterize the functional properties of the associated weighted Lorentz spaces.


\medskip


\newpage
\section{Introduction}

Recently, some authors (see \cite{Ba}, \cite{BPSo}, and \cite{BPSt})  have considered
multidimensional analogs of classical inequalities for monotone functions: Hardy's inequality,
Chebyshev's inequality, embeddings for weighted Lorentz spaces, etc. (see, e.g.,  \cite{AM},
\cite{Sa},  \cite{St},
\cite{CS}).  We recall that the main interest in studying these results on monotone functions comes
from the fact that the  spaces, where the  estimates hold, are rearrangement invariant function
spaces (see \cite{BS}), and hence the functions that show up in the inequalities are the nonincreasing
rearrangements of general measurable functions (which are essentially all monotone functions on
$\RR_+$).  This observation is fundamental to understand our main purpose: we want to find the
natural definition for a multidimensional rearrangement in such a way that what we get is a general
decreasing function on $\RR^n_+:=\RR_+\times\cdots\times\RR_+$. Our approach is very
geometrical: we look for a measure preserving  transformation taking (all) sets in $\RR^n$ to (all)
decreasing sets in
$\RR^n_+$, and such that it is monotone, and leaves fixed the sets that are already decreasing (see
Definition~\ref{set}). Once we know how to rearrange sets, we can define the multidimensional
rearrangement of a function by using the \lq\lq Layer-cake formula\rq\rq, which recovers a function
by means of its level sets (see Definition~\ref{rearange}).

\medskip

This new definition opens up the possibility of studying whether the  properties of the classical
rearrangement hold true in  the multidimensional setting (see Corollary~\ref{harlit} for an example
which shows that  the resonant property fails). In Section~\ref{twodimdr} we develop the main ideas
of the new rearrangement from a measure theoretical point of view (Propositions~\ref{propset},
\ref{propre},
\ref{simpfun}, and Theorem~\ref{betharli}), establish the relationship with the classical
rearrangement and show that it agrees with the so called multivariate rearrangement
(Corollary~\ref{equirear} and Theorem~\ref{iterrea}). In Section~\ref{lorspawd} we introduce the
weighted Lorentz spaces associated to the multidimensional rearrangement, we find their relationship
with the Lebesgue and the rearrangement invariant spaces (Theorem~\ref{lorleb} and
Propositions~\ref{wconst},
\ref{exwlp}), prove the different embeddings in the whole range of indices
(Proposition~\ref{embddg}), and characterize functional properties like quasinormability
(Theorem~\ref{quasiw}) and the weights  which give rise to a norm (Theorems~\ref{indexp} and
\ref{normeq}).

\medskip

Most of the notations we are going to use are standard as, for example, defined in \cite{BS}:
$\lambda_f$ is the distribution function of $f$, the nonincreasing rearrangement of $f$  is
denoted $f^*$,
$h\downarrow$ means that $h$ is decreasing,  etc. A weight
$w$ is a locally integrable positive function (either on $\RR^n$ or $\RR^n_+$, depending on the
context), and if $E$ is a set, $w(E)=\int_Ew$. As usual, $|E|$ denotes the Lebesgue measure of $E$.
Two  positive quantities $A$ and $B$, are said to be equivalent
($A\approx B$) if there exists a constant $C>1$ (independent of the essential parameters 
defining $A$ and
$B$, and not the same at different occasions) such that $C^{-1}A\le B\le CA$. Also, all sets that we
are going to consider are always Lebesgue measurable sets.

\bigbreak
\section{Two-dimensional decreasing rearrangement}\label{twodimdr}

For simplicity, we are going to reduce our definitions to the two-dimensional case, although there are
natural extensions to higher dimensions too. Our approach is to give a geometric definition of the
rearrangement of a  measurable set (so that we get a general decreasing set in ${\RR}^2_+$),
and extend it to also rearrange functions, by looking at the level sets and the use of the Layer-cake
formula (\cite{LL}). We will show in Theorem~\ref{iterrea} that this definition agrees with the so
called multivariate rearrangement (see \cite{Bl}).

\begin{defi}\label{def1}  \rm We say that a set $D\subset{\RR}^2_+$ is
decreasing (and write $D\in\Delta_d$) if the function $\chi_D$ is
decreasing in each variable. 
\end{defi}

\begin{defi}\label{set}  \rm
 Let $E$ be a subset of ${{\RR}^2}$ and
$\varphi_E(x)=|\{y\in {\RR}:(x,y)\in E\}|$,
 $x\in {\RR}$. Let the function $\varphi_E^{*}$, defined by
$$
\varphi_E^{*}(s)=\inf\{\lambda:|\{x\in {\RR}
: \varphi_E(x) >\lambda\}|\le s\},\,(s\ge 0)
$$
be the usual decreasing rearrangement of $\varphi_E$ (see \cite{BS}). Then, the two-dimensional
decreasing rearrangement of the set $E$  is the set
$$
E^{*}=\{(s,t)\in {{\RR}^2_+}:0<t<\varphi_E^{*}(s)\}.
$$
\end{defi}

\begin{defi}\label{rearange} {\bf (Layer-cake formula \cite{LL}).}  \rm
 The two-dimensional decreasing rearrangement $f^{*}_2$ of a
function $f$ on ${{\RR}^2}$ is given by
$$
f^{*}_2(x)=\int_0^{\infty}\chi_{\{| f|>t\}^{*}}(x) dt,\qquad x\in\RR^2_+.
$$
\end{defi}

We give now some elementary properties for this new rearrangement definition.  

\begin{prop}\label{propset}   
Let $E$ and $F$  be two   subsets of ${{\RR}^2}$. 
Then, 

\medskip
\noindent
a) $| E| =| E^{*}|$, and $E^*\subset F^*$, if $E\subset F$.

\medskip
\noindent
b) $E= E^{*}$, if and only if $E$ is a decreasing set of  ${{\RR}_+^2}$.

\medskip
\noindent
c) $\fst=\chi_{F^{*}}$, if and only if 
$f=\chi_E$, and $E^{*}=F^{*}$. In
particular, $\left({\chi_E}\right)^{*}_2=\chi_{E^{*}}.$

\medskip
\noindent
d) If $E\cap F=\emptyset$ then  $|(E\cup F)^{*}\setminus E^{*}|=|F|$.
\end{prop}

{\bf Proof.} 
a) We have 
$$
| E| = \int_{-\infty}^{\infty}\varphi_E(x)
dx=\int_{0}^{\infty}{\varphi}_E^{*}(x) dx=| E^{*}|.
$$
The second part is trivial since $\varphi_E\le\varphi_F$.

\medskip
\noindent
b) If $E$ is a decreasing set, then there exists $r>0$ such that
$$
E=\{(x,y)\in {{\RR}^2}:0<x<r,\,0<y<{\varphi_E}(x)\}.
$$
Since $\varphi_E$ is decreasing, then  $E= E^{*}$. The converse implication is trivial.

\medskip
\noindent
c) It yields that
$$
\left({\chi_E}\right)^{*}_2 (x)= \int_0^{\infty}
\chi_{\{\chi_E>t\}^{*}}(x) dt
 =\int_0^{1}\chi_{E^{*}}(x) dt=\chi_{E^{*}}(x). 
$$
Conversely, suppose $\fst=\chi_{F^*}$:
\begin{enumerate}
\item[--]
If $x\notin F^*$, then $x\notin\{f>t\}^*$ and hence $\{f>t\}^*\subset F^*$, for all $t>0$.

\item[--]
If $x\in F^*$, then $x\in\{f>t\}^*$, $0<t<1$, and $x\notin\{f>t\}^*$, $1<t$.
\end{enumerate}

Therefore, $\{f>t\}^*= F^*$, if $0<t<1$, and $\{f>t\}^*=\emptyset$, if
$1<t$. Thus, $t<f(x)\le 1$, if $f(x)\neq0$, for every $0<t<1$, and hence there exists a set $E$
such that $f=\chi_E$ and $E^*=F^*$.

\medskip
\noindent
Property d) follows easily from a).   $\qquad\Box$
\medskip

The following results gives more information on the level sets of $f$ and $\fst$.

\begin{lemma}\label{incl}
If $f$ is a measurable function on ${{\RR}^2}$ and $t>0$, then

$$\{f^{*}_2>t\}\subseteq \{| f|>t\}^{*} \subseteq \{f^{*}_2\ge t\}.$$
\end{lemma}

{\bf Proof.}  By definition,
$$ 
f^{*}_2(x)>t\Longleftrightarrow \int_0^{\infty}\chi_{\{| f|>s\}^{*}}(x)
ds>t,\,\,(x=(x_1,x_2)).
$$
But,
$$
\chi_{\{| f|>s\}^{*}}(x)
=\left\{ \begin{array}{ll}
1 & \mbox{if ${\varphi}_s^{*}(x_1)>x_2$}\\
  & \mbox{}\\
0 & \mbox{if ${\varphi}_s^{*}(x_1)\le x_2,$}
\end{array}
\right.
$$
where ${\varphi}_s(a)=|\{b:| f(a,b)| >s\}|$.
Thus $f_2^{*}(x)>t \Longleftrightarrow
|\{s:{\varphi}_s^{*}(x_1)>x_2\}|>t$. 
Observe that if
$
s<s^{\prime}$, then $   {\varphi}_{s^{\prime}}^{*}(x_1)\le
{\varphi}_s^{*}(x_1),$
and hence $\{s:{\varphi}_s^{*}(x_1)>x_2\}$ is an interval of the form
$(0,s)$ or $(0,s]$. 
Hence,
\begin{eqnarray}
|(0,s)| >t  & \Longrightarrow & s>t  \Longrightarrow
{\varphi}_t^{*}(x_1)>x_2\nonumber\\
   & \Longrightarrow & (x_1,x_2)=x\in \{ | f| >t\}^{*}.\nonumber
\end{eqnarray}
Conversely, if $ x\in \{| f |>t\}^{*}$, then
$\varphi_t^{*}(x_1)>x_2,\, x=(x_1,x_2),$
and hence $|\{s:{\varphi}_s^{*}(x_1)>x_2\}|\ge t$, which implies $f^{*}_2(x)\ge
t.\qquad\Box$

\medskip

\begin{lemma}\label{car}
Let $f$ and $g$ be two measurable functions on ${{\RR}^2}$ and $t>0$.Then
$$
\chi_{\{| f+g| > t\}^{*}} (x+y) \le \chi_ {\{| f| > t/2
\}^{*}}(x)
+\chi_ {\{| g| > t/2 \}^{*}}(y),
$$
$x=(x_1,x_2),\,\,y=(y_1,y_2)\in {{\RR}^2}$.
\end{lemma}

{\bf Proof.}
Let
\begin{eqnarray*}
\varphi_{f,t}(a)&:=&|\{b\in {\RR}:| f(a,b)|>t\}|\\
\varphi_{g,t}(a)&:=&|\{b\in {\RR}:| g(a,b)|>t\}|\\
\varphi_{f+g,t}(a)&:=&|\{b\in {\RR}:| (f+g)(a,b)|>t\}|.
\end{eqnarray*}
We know that
$$
\varphi_{f+g,t}(a) \le \varphi_{f,t/2}(a)+\varphi_{g,t/2}(a).
$$
Also, if $x\notin\{ | f |>t/2\}^{*}$, then $
\varphi_{f,t/2}^{*}(x_1)<x_2$ and similarly, if
$y\notin\{ | g |>t/2\}^{*}$, then $
\varphi_{g,t/2}^{*}(y_1)<y_2$. 
Therefore
\begin{eqnarray*}
{\varphi}_{f+g,t}^{*}(x_1+y_1) & \le &
(\varphi_{f,t/2}+\varphi_{g,t/2})^{*}(x_1+y_1)\nonumber\\
 & \le & \varphi_{f,t/2}^{*}(x_1)+\varphi_{g,t/2}^{*}(y_1)\nonumber\\
 & < & x_2+y_2,
\end{eqnarray*}
which means exactly that
$ x+y\notin\{ | f+g|> t\}^{*}$. This completes the proof
of the lemma.$\qquad\Box$

\begin{prop}\label{propre}
Suppose $f$, $g$, and $f_n$, ($n=1,2,\ldots$) are measurable functions on
${{\RR}^2}$ and let $c\in \CC$. 
 Then the two-dimensional decreasing rearrangement $f^{*}_2$ is a
nonnegative function  on ${{\RR}^2_+}$, decreasing in each variable. Furthermore,

\medskip
\noindent
a) $ | g|\le | f| -\mbox{a.e.}\Longrightarrow g^{*}_2\le
f^{*}_2;$

\medskip
\noindent
b)
$ (cf)^{*}_2=| c| f^{*}_2; $

\medskip
\noindent
c) if $f$ is decreasing in each variable, then $f^{*}_2=f;$

\medskip
\noindent
d) $ (f+g)^{*}_2(x+y)\le 2\left(f^{*}_2(x)+g^{*}_2(y)\right);$

\medskip
\noindent
e)
$ | f| \le \liminf_{n\rightarrow \infty} | f_n|
\Longrightarrow f^{*}_2\le
 \liminf_{n\rightarrow \infty}\left(f_n\right)_2^{*} $, and, 
in particular, if 
$| f_n| \uparrow | f| \Longrightarrow
\left(f_n\right)_2^{*}\uparrow f^{*}_2$;

\medskip
\noindent
f)
$  \left(f^{*}_2(x)\right)^p=
\left(f^p(x)\right)^{*}_2,\,\,(0<p<\infty) $;

\medskip
\noindent
g) if $f$ is a symmetric function (i.e. $f(x_1,x_2)=f(x_2,x_1)$), then
$f^{*}_2$ is symmetric.
\end{prop}

{\bf Proof.} That $f^{*}_2$ is nonnegative and decreasing follows from
Definition~\ref{rearange} and the fact
 that the characteristic function of a decreasing set is a decreasing
function.

\medskip
\noindent
a) By Definition~\ref{set}, it follows that
$$
\{| g| >t\}\subset \{| f| >t\} \Longrightarrow \{| g|
>t\}^{*}\subset \{| f| >t\}^{*}.
$$
Thus
$
\chi_{\{| g| >t\}^{*}} \le \chi_{\{| f| >t\}^{*}}
$
and
$
g^{*}_2\le f^{*}_2.
$

\medskip
\noindent
b) Trivial.

\medskip
\noindent
c) If $f$ is a decreasing function in each variable, then the level set $ \{| f| >t\}$ is a
decreasing set (see also \cite{BPSo})
 and c.f. Proposition~\ref{propset}
$$
 \{| f| >t\}^{*}= \{| f| >t\}.
$$
We get the desired equality by using Definition~\ref{rearange}.

\medskip
\noindent
d) By Lemma~\ref{car} and   b) of this proposition we have
\begin{eqnarray}
(f+g)^{*}_2(x+y) & = &\int_0^\infty\chi_{\{| f+g |
>t\}^{*}}(x+y) dt\nonumber\\
& \le & \int_0^\infty\chi_{\{| f| >t/2\}^{*}}(x) dt+\int_0^\infty
\chi_{\{| g| >t/2\}^{*}}(y)dt\nonumber\\
 & = & 2\left(f^{*}_2(x)+g^{*}_2(y)\right).\nonumber
\end{eqnarray}

\medskip
\noindent
e) Let
$$
E^t:=\{(x,y):| f(x,y)|>t\}
$$
and
$$
E^t_n:=\{(x,y):| f_n(x,y)|>t\}.
$$
Set $f_x(y):=f(x,y)$ and
$$
\varphi_{f,t}(x)=|\{y:| f(x,y)| > t\}|= \lambda_{f_x}(t),
$$
where $\lambda_{f_x}$ is the usual distribution function (see \cite{BS}). Then
\begin{eqnarray*}
| f| \le \liminf_{n\rightarrow \infty} | f_n| & \Longrightarrow &
 | f_x| \le \liminf_{n\rightarrow \infty} | f_{x,n}|\,\,
\mbox{a.e.}\nonumber\\
& \Longrightarrow & \lambda_{f_x} \le \liminf_{n\rightarrow \infty}
\lambda_{f_{x,n}}\nonumber\\
& \Longrightarrow &\varphi_{f,t} \le \liminf_{n\rightarrow \infty}
\varphi_{f_n,t}, \mbox{a.e.},\,\forall t>0\nonumber\\
& \Longrightarrow & {\varphi}^{*}_{f,t} \le \liminf_{n\rightarrow
\infty} {\varphi}^{*}_{f_n,t}, \mbox{a.e.},
\,\forall t>0\nonumber\\
& \Longrightarrow & \chi_{\left(E^t\right)^{*}} \le
\liminf_{n\rightarrow \infty}
\chi_{\left(E_n^t\right)^{*}}\nonumber\\
& \Longrightarrow & f^{*}_2\le \liminf_{n\rightarrow
\infty}\left(f_n\right)_2^{*}.
\end{eqnarray*}
The second part is an immediate consequence of the first.

\medskip
\noindent
f) We have
 $$ \left(f^{*}_2(x)\right)^p=\left(\int_0^{\infty}\chi_{\{|
f|>t\}^{*}}(x)dt\right)^p$$
\begin{eqnarray}
\left(f^p\right)^{*}_2(x) & = & \int_0^{\infty} \chi_{\{|
f^p|>t\}^{*}}(x)dt\nonumber\\
& = & p\int_0^{\infty} \chi_{\{| f|>t\}^{*}}(x)t^{p-1}dt.\nonumber
\end{eqnarray}
In view of  Lemma~\ref{incl}  we have
$$
\chi_{\{| f|>t\}^{*}}\ge \chi_{\{f^{*}_2\ge t\}}
$$
and, hence,
\begin{eqnarray}
\left(f^p\right)_2^{*} (x) & \ge &
p\int_0^{\infty}t^{p-1}\chi_{\{f^{*}_2\ge t\}}(x) dt \nonumber\\
& = &
p\int_0^{f^{*}_2(x)}t^{p-1}dt=\left(f^{*}_2(x)\right)^{p}.\nonumber
\end{eqnarray}
On the other hand, if we take $0<r<1$, then by Lemma~\ref{incl} we have
$$
\chi_{\{| f|>t\}^{*}_2} \le \chi_{\{f^{*}_2\ge r t\}},
$$
and, hence,
\begin{eqnarray}
\left(f^p\right)_2^{*} (x) & \le &
p\int_0^{\infty}t^{p-1}\chi_{\{f^{*}_2\ge r t\}}(x) dt \nonumber\\
& = &
p\int_0^{f^{*}_2(x)/r}t^{p-1}dt=\left(\frac{f^{*}_2(x)}{r}\right)^{p}.\nonumber
\end{eqnarray}
Since this is true for all $0<r<1$, we get
$$
f^{*}_2(x)\ge \left(f^p\right)_2^{*} (x) \ge
f^{*}_2(x).
$$

\medskip
\noindent
g) This is just an observation which follows immediately by using the definition of $f^{*}_2$.
$\qquad\Box$

\medskip

The following proposition will be very useful for proving our main results, since it will allow us to
consider the special and easier case of simple functions. 

\begin{prop}\label{simpfun}
If $f$ is a measurable function on ${{\RR}^2}$, then there exists a sequence $(s_n)_n$ of simple
measurable functions such that:

\medskip
\noindent
a) $0\le (s_1)^{*}_2 \le \ldots \le (s_n)^{*}_2 \le f^{*}_2$,

\medskip
\noindent
b) $(s_n)^{*}_2\longrightarrow f^{*}_2$ as $n\rightarrow \infty$
a.e.
\end{prop}

{\bf Proof.} The existence of the sequence is standard, and the rest is just
a consequence of Proposition~\ref{propre} a) and e), and the following remark: If
$s(x)=\sum_{j=1}^na_j\chi_{E_j}$, with $a_1>a_2>\cdots>a_n>0$ and $E_j\cap
E_i=\emptyset,\,i\ne j$, then 
$$
s^*_2(x)=\sum_{j=1}^na_j\chi_{F^*_j\setminus F^*_{j-1}}(x),
$$
where $F_j=\cup_{k=1}^jE_k$, and $F_0=\emptyset$. Observe that from Proposition~\ref{propset}
we have that $|F^*_j\setminus F^*_{j-1}|=|E_j|$.  $\qquad\Box$
\medskip

As a corollary, we can obtain several properties relating our two-dimensional rearrangement and the
classical one. In particular, we see that the new rearrangement is finer and gives more information
than the other.

\begin{cor}\label{equirear}
Let $f$ and $g$ be two measurable functions in $\RR^2$.  

\medskip
\noindent
a) If $f^*_2=g^*_2$, then $f^*=g^*$, and the converse is not true in general.

\medskip
\noindent
b) $(f^*_2)^*=f^*.$

\end{cor}

{\bf Proof.} To prove a) we observe that if $f^*_2=g^*_2$, then 
$$
\int_0^{\infty}\chi_{\{f>t\}^*}(x)\,dt=\int_0^{\infty}\chi_{\{g>t\}^*}(x)\,dt,
$$
and hence $\{f>t\}^*=\{g>t\}^*$. Using now Proposition~\ref{propset} a), we get that
$|\{f>t\}|=|\{g>t\}|$ which shows that $f^*=g^*$. To see that the converse does not hold, consider
the decreasing sets $A=(0,1)\times(0,2)$, $B=(0,2)\times(0,1)$ and the functions $f=\chi_A$ and
$g=\chi_B$. Then, $f^*=g^*=\chi_{(0,2)}$ but $f^*_2=f\neq g=g^*_2$.

\medskip
The proof of b) follows immediately by checking what happens for simple functions and using
Proposition~\ref{simpfun}.  We observe that from b) we can  also give an alternative proof of
a)  $\qquad\Box$

\medskip

We consider next integral inequalities, for the two-dimensional rearrangement, related to the
Hardy--Littlewood inequality (see \cite{BS}). Again we observe that what we obtain is a better
estimate. We begin with an elementary but useful lemma.

\medskip

\begin{lemma}\label{haliprop}
Let $g$ be a  nonnegative simple function on ${{\RR}^2}$  and let $E$ be an arbitrary  set
of
${{\RR}^2}$. Then
$$
\int_E g(x)\, dx\le \int_{E^{*}} g^{*}_2(x)\, dx.
$$
\end{lemma}

{\bf Proof.} Let
$$
g(x)=\sum_{j=1}^n a_j\chi_{E_j}(x),
$$
where $a_1>a_2>\cdots>a_n>0$,  $a_{n+1}=0$, and $E_j\subset {{\RR}^2}$ are
 of finite measure such that
 $E_j\cap E_i=\emptyset,\,i\ne j$. 
Another representation of $g$ is
$$
g(x)=\sum_{j=1}^n b_j\chi_{F_j}(x),
$$
where $b_j>0$, $b_j=a_j-a_{j+1}$,   and $F_j={\cup}_{i=1}^j E_i$. Then, 
\begin{eqnarray}
g^{*}_2(x) & = & \int_0^{\infty} \chi_{\{g> t\}^{*}}(x) dt
\nonumber\\
& = & \int_{a_2}^{a_1}\chi_{{E_1}^{*}}(x)
dt+\int_{a_3}^{a_2}\chi_{{(E_1\cup E_2})^{*}}(x) dt
+\cdots \int_0^{a_n}\chi_{{(E_1\cup\ldots \cup E_n)}^{*}}(x) dt
\nonumber\\
& = & \chi_{{E_1}^{*}}(x) (a_1-a_2)+\chi_{{(E_1\cup E_2)}^{*}}(x)
(a_2-a_3)+\cdots+
\chi_{{(E_1\cup\ldots \cup E_n)}^{*}}(x)a_n\nonumber\\
& = & \sum_{j=1}^n b_j\chi_{{F_j}^{*}}(x).\label{reofsim}
\end{eqnarray}

\noindent
Thus, since  $(F_j\cap E)^{*}\subset {F_j}^{*}\cap
E^{*}$,  we have that
\begin{eqnarray}
\int_E g(x) dx & = & \sum_{j=1}^n b_j\int_E\chi_{F_j}(x) dx=\sum_{j=1}^n
b_j| F_j\cap E|\nonumber\\
& = & \sum_{j=1}^n b_j| (F_j\cap E)^{*} | =\sum_{j=1}^n b_j
\int_{(F_j\cap E)^{*} }dx\nonumber\\
& \le & \sum_{j=1}^n b_j \int_{ F_j^{*}\cap E^{*}}dx =
\int_{E^{*}} g^{*}_2(x) dx. \qquad\Box\nonumber
\end{eqnarray}

\medskip

\begin{thr}\label{betharli}
If $f$ and $g$ are measurable functions on ${{\RR}^2}$, then
$$
\int_{{\RR}^2}| f(x) g(x)|\, dx\le \int_{{\RR}^2_+}
f^{*}_2(x) g^{*}_2(x)\, dx\le \int_0^{\infty}
f^{*}(t) g^{*}(t)\, dt.
$$
\end{thr}

{\bf Proof.} It is enough to prove the statement for $f$ and $g$ nonnegative. By
Proposition~\ref{propre} e) and in
view of the monotone convergence theorem there is no loss of generality in
assuming $f$ and $g$ to be simple functions. Let
$$
f(x)=\sum_{j=1}^n a_j\chi_{E_j}(x),
$$
where $E_1\subset E_2\subset\ldots\subset\ldots E_n\subset{{\RR}^2}$,
are  
of finite measure, and $a_j>0$. Then, by Lemma~\ref{haliprop}, we have that 
\begin{eqnarray}
\int_{{\RR}^2} f(x)g(x)dx & = & \sum_{j=1}^n a_j\int_{E_j} g(x) dx
\leq \sum_{j=1}^n a_j\int_{E_j^{*}}
 g^{*}_2(x)dx \nonumber\\
& = &\int_{{\RR}^2_+} \sum_{j=1}^n a_j \chi_{E_j^{*}}(x)
g^{*}_2(x)dx = \int_{{\RR}_+^2} f^{*}_2(x)
 g^{*}_2(x) dx .  \nonumber
\end{eqnarray}

\noindent
The second inequality follows from Corollary~\ref{equirear} b).  $\qquad\Box$

\medskip

\begin{cor}\label{harlit}

If $f$ is a nonnegative measurable function on ${{\RR}^2}$, and $D$ is a decreasing set, then

$$
\sup_{E^*=D}\int_Ef(x)\,dx\le\int_Df^*_2(x)\,dx\le\int_0^{|D|}f^*(t)\,dt,
$$
and both inequalities can hold strictly for some $f$ and $D$.

\end{cor}

{\bf Proof.} That the inequalities hold is a consequence of Theorem~\ref{betharli}, applied with
$g=\chi_E$. To show that the first inequality can be strict, consider the sets $A=(3,4)\times(0,1)$,
$B=(4,6)\times (0,2)$, $D=(0,1)\times(0,2)$, and the function $f(x)=2\chi_A(x)+\chi_B(x)$.
Then, it is easy to see that for every set $E$ such that $E^*_2=D$, we have 
$$
\int_Ef(x)\,dx\le 2<3=\int_Df^*_2(x)\,dx.
$$

\noindent
For the second inequality, consider $D_{\varepsilon}=(0,\varepsilon)\times(0,1/\varepsilon)$ and
$f$ as before. Then, 
$$
\int_0^{|D_{\varepsilon}|}f^*(t)\,dt=2,\qquad\hbox{for every $\varepsilon>0$,}
$$
but
$$
\lim_{\varepsilon\to 0}{\int_{D_{\varepsilon}} f^*_2(x)\,dx=0. }\qquad\Box
$$

\medskip

As we have mentioned in the introduction, our definition of the two-dimensional rearrangement is based
on a geometric approach: we first look at the rearrangement of the level sets of the function, and then
we recover the rearrangement of the function by summing up all these level sets (Layer-cake
formula). In the next theorem, we are going to prove a direct way of calculating the two-dimensional
rearrangement as an iterative procedure with respect to the usual rearrangement in each
variable (see \cite{Bl} for some related work). 

\medskip

In order to clarify the notation used in the proof, given a function $f(x,y)$ defined on $\RR^2$, we
write $R_t(x)=(f_x)^{*y}(t)$, where $f_x(y)=f(x,y)$ and $t>0$ (i.e., $R_t$ is the usual
rearrangement of the function $f_x$, with respect to the variable $y$). Similarly, we set $\tilde
f(s,t)=(R_t)^{*x}(s)$, $s,t>0$. It is very easy to show that, in general, we do not get the same
function if we first rearrange with respect to $x$ and then with respect to $y$.

\medskip

\begin{thr}\label{iterrea}
If $f$ is a measurable function on ${{\RR}^2}$, then $f^*_2(s,t)=\tilde f(s,t)$, $\forall\ s,t>0$.

\end{thr}

{\bf Proof.} Using Proposition~\ref{simpfun}, it suffices to consider $f$ to be a simple function.
Hence, let $f(x,y)=\sum_{j=1}^na_j\chi_{E_j}(x,y)$, with $a_1>a_2>\cdots>a_n$, $E_j\cap
E_k=\emptyset$, $j\neq k$. Set $F_k=\cup_{j=1}^kE_j$, $F_0=\emptyset$, so that 
$$
f^*_2(s,t)=\sum_{j=1}^na_j\chi_{F^*_j\setminus F^*_{j-1}}(s,t).
$$
Recall that $\varphi_E(x)=|\{y:(x,y)\in E\}|$ and $E^*=\{(s,t):0<t<\varphi^*_E(s)\}$. Hence,
$$
\chi_{E^*}(s,t)=\chi_{(0,\varphi^*_E(s))}(t)=\chi_{(0,\lambda_{\varphi_E}(t))}(s).
$$
Thus,
$$
\chi_{F^*_j\setminus F^*_{j-1}}(s,t)=\chi_{F^*_j}(s,t)-\chi_{
F^*_{j-1}}(s,t)=\chi_{[\lambda_{\varphi_{F_{j-1}}}(t),\lambda_{\varphi_{F_j}}(t))}(s),
$$
which gives
\begin{eqnarray}\label{eqfst}
f^*_2(s,t)=\sum_{j=1}^na_j\chi_{[\lambda_{\varphi_{F_{j-1}}}(t),\lambda_{\varphi_{F_j}}(t))}(s).
\end{eqnarray}

On the other hand, since
$$
f_x(y)=\sum_{j=1}^na_j\chi_{E_j(x)}(y),
$$
where $E(x)=\{y:(x,y)\in E\}$, we have that

\begin{eqnarray*}
R_t(x)&=&(f_x)^{*y}(t)=\sum_{j=1}^na_j\chi_{[|F_{j-1}(x)|,|F_j(x)|)}(t)\\
&=&\sum_{j=1}^na_j\chi_{[\varphi_{F_{j-1}}(x),\varphi_{F_{j}}(x))}(t)
=\sum_{j=1}^na_j\chi_{H_j(t)}(x),
\end{eqnarray*}
where $H_j(t)=\{y:\varphi_{F_{j-1}}(y)\le t<\varphi_{F_{j}}(y)\}$. Therefore,

\begin{eqnarray}\label{eqfti}
\tilde f(s,t)=(R_t)^{*x}(s)=\sum_{j=1}^na_j\chi_{[|G_{j-1}(t),|G_j(t)|)}(s),
\end{eqnarray}
where $G_j(t)=\cup_{k=1}^jH_k(t)$, $G_0(t)=\emptyset$. Thus looking at (\ref{eqfst}) and
(\ref{eqfti}) it suffices to proving that
$$
|G_j(t)|=\lambda_{\varphi_{F_j}}(t).
$$
But, in fact
\begin{eqnarray*}
|G_j(t)|&=&\sum_{k=1}^j|H_k(t)|=\sum_{k=1}^j|\{y:\varphi_{F_{k-1}}(y)\le
t<\varphi_{F_{k}}(y)\}|\\ 
&=&|\{y:t<\varphi_{F_j}(y)\}|=\lambda_{\varphi_{F_j}}(t),
\end{eqnarray*}
and the proof is complete. $\qquad\Box$
\medskip

\begin{cor}
If $g$ and $h$ are two measurable functions on $\RR$, and $f(x,y)=g(x)h(y)$, then
$f^*_2(s,t)=g^*(s)h^*(t).$
\end{cor}

\medskip

Another application of Theorem~\ref{iterrea} is that the inequality proved in Theorem~\ref{propre}
d) can be improved to obtain the classical subadditivity condition:  $ (f+g)^{*}_2(x+y)\le
f^{*}_2(x)+g^{*}_2(y)$ (we leave the details to the interested reader).

\bigbreak

\section{A new multidimensional Lorentz space}\label{lorspawd}

In this section we prove some properties of a new type of space, defined
using the two-dimensional decreasing rearrangement. Recall the definition of the classical Lorentz
space: If $v$ is a weight in $\RR_+$ and $0<p<\infty$,
$$
\Lambda^p(v)=\bigg\{f:\RR^n\to\CC:\Vert
f\Vert_{\Lambda^p(v)}:=\bigg(\int_0^{\infty}(f^*(t))^pv(t)\,dt\bigg)^{1/p}<\infty\bigg\}.
$$ 
We now say that a measurable function
$f$ on
${{\RR}^2}$ belongs to the (multidimensional) Lorentz space
${\Lambda}^p_2(w)$, provided $\|f\|_{{\Lambda}^p_2(w)}$,
defined by
\begin{equation}
\| f\|_{{\Lambda}^p_2(w)}:=\left(\int_{{\RR}_+^2}
\left(f^{*}_2(x)\right)^p w(x) dx\right)^{1/p},
\end{equation}
is finite.
Here $w$ is a nonnegative, locally integrable function on ${{\RR}_+^2}$, not identically $0$.

The next result gives an alternative description of the $L^p_{{\RR}^2}$ norm in terms
of the two-dimensional decreasing rearrangement,
 i.e., the spaces defined above generalize naturally the Lebesgue spaces.

\begin{thr}\label{lorleb}
 If $0<p<\infty$, then
$
{\Lambda}^p_2(1)=L^p_{{\RR}^2}.
$
\end{thr}

{\bf Proof.} By Fubini's theorem and Proposition~\ref{propre} f) we have
\begin{eqnarray}
\int_{{\RR}^2}|f(x)|^p dx & = &
\int_0^{\infty}\int_{\{|f|^p>t\}}dx dt=\int_0^{\infty}\int_{\{|f|^p>t\}
^{*}}dx dt\nonumber\\
& = & \int_{{\RR}_+^2}\int_0^{\infty}\chi_{\{|f|^p>t\}^{*}}(x) dt dx
 = \int_{{\RR}_+^2} \left(f^p\right)^{*}_2 (x) dx \nonumber\\
& = & \int_{{\RR}_+^2} \left(f^{*}_2(x)\right)^p dx . \qquad\Box\nonumber
\end{eqnarray}

\medskip

We are interested in studying functional properties of the spaces ${\Lambda}^p_2(w)$ and their
relationship with the classical rearrangement invariant spaces (see \cite{BS}). The following results
show that these two kinds of spaces only agree in very particular cases:

\begin{prop}\label{wconst}
If $\Vert\cdot\Vert_{{\Lambda}^p_2(w)}$ is a rearrangement invariant norm, then $w$ is constant,
and hence ${\Lambda}^p_2(w)=L^p_{{\RR}^2}.$
\end{prop}

{\bf Proof.} Fix $(x,y)\in\RR^2_+$, $0<\varepsilon<\min(x,y)$,  and define
$R=(0,x)\times(0,y)$, 
$P_{\varepsilon}=(x-\varepsilon,x)\times (y-\varepsilon,y)$,
$Q_{\varepsilon}=(x,x+\varepsilon)\times(0,\varepsilon)$, and
$A_{\varepsilon}=(R\setminus P_{\varepsilon})\cup Q_{\varepsilon}$. Then
$|R|=|A_{\varepsilon}|$, and hence
$\Vert\chi_R\Vert_{{\Lambda}^p_2(w)}=\Vert\chi_{A_{\varepsilon}}\Vert_{{\Lambda}^p_2(w)}$,
which gives
$$
\int_{P_{\varepsilon}}w(x)\,dx=\int_{Q_{\varepsilon}}w(x)\,dx.
$$
Now, letting $\varepsilon\to0$, using the Lebesgue differentiation theorem, and a symmetric
argument changing $x$ and $y$, we obtain that
$w$ is constant.  $\qquad\Box$

\medskip

In a similar way, one can prove the following:

\begin{prop}\label{exwlp}
There exists a weight $v$ in $\RR^2$ such that $\Lambda^p_2(w)=L^p_{\RR^2}(v)$ if and only if
$\Lambda^p_2(w)=L^p_{\RR^2}$.
\end{prop}

\medskip

It is very easy to see that embedding results for the spaces $\Lambda^p_2(w)$ are equivalent to
embeddings  for the cone of decreasing functions on $L^p_{\RR^2_+}$, which have been completely
characterized in all cases (see \cite{BPSo} and \cite{BPSt}). The result reads as follows:

\begin{prop}\label{embddg} Let $0<p_1,p_2<\infty$ and $w_1,w_2$ be two weights in $\RR^2_+$.

\medskip
\noindent
a) If $p_1\le p_2$, then $\Lambda^{p_1}_2(w_1)\subset\Lambda^{p_2}_2(w_2)$, if and only if,
$$
\sup_{D\in\Delta_d}{w_2(D)^{1/p_2}\over w_1(D)^{1/p_1}}<\infty.
$$

\medskip
\noindent
b) If $p_1>p_2$, then $\Lambda^{p_1}_2(w_1)\subset\Lambda^{p_2}_2(w_2)$, if and only if,
$$
\sup_{0\le
h\downarrow}\int_0^{\infty}w_1(D_{h,t})^{-r/p_1}d(-w_2(D_{h,t})^{r/p_2})<\infty,
$$
where $D_{h,t}=\{x\in\RR^2_+:h(x)>t\}$, and $1/r=1/p_2-1/p_1$.
\end{prop}

\medskip

The characterization of the quasinormability, in the case of the classical Lorentz spaces,  was proved in
\cite{CS} to be equivalent to a doubling condition on the weight (the $\Delta_2$-condition).  We show
that a similar result holds for the two-dimensional rearrangement.

First we note that the spaces $\Lambda_2^p(w)$, $0<p<\infty$, have the following (quasi)norm
properties:
\begin{equation}\label{hom}
\Vert cf\Vert_{\Lambda_2^p(w)}=|c|\Vert f\Vert_{\Lambda_2^p(w)},
\end{equation}
(see Proposition~\ref{propre} b)), and if $w$ is strictly positive (which we assume in the sequel)
\begin{equation}\label{can}
\Vert f\Vert_{\Lambda_2^p(w)}=0\Longleftrightarrow f=0 \hbox{ a.e.}
\end{equation}
Thus, in order to investigate if $\Vert\cdot\Vert_{\Lambda_2^p(w)}$ is a norm (quasi-norm) we
only have to check that the triangle (quasi-triangle) inequality holds.
\medskip

\begin{thr}\label{quasiw}
Let $0<p<\infty$. Then, 
$\| \cdot\|_{{\Lambda}^p_2(w)}$ is a quasinorm if and only if there exists a
constant $C>0$ such that
\begin{equation}\label{quasi}
\int_D w(2x) dx \le C \int_D w(x) dx, 
\end{equation}
for all decreasing sets  $D\subset {\RR}_+^2$.  
Moreover, with this quasinorm, ${{\Lambda}^p_2(w)}$ becomes a complete
quasinormed space.
\end{thr}

{\bf Proof.} For sufficiency we use Proposition~\ref{propre} d), Theorem
2.2 d) in \cite{BPSo}, with $p=q$, and we get:

\begin{eqnarray}
\| f+g\|^p_{{\Lambda}^p_2(w)} & = & \int_{{\RR}_+^2}\left(\left(f+g\right)_2^{*}(x)\right)^p w(x)
dx\nonumber\\
& \le & C\int_{{\RR}_+^2}
\left(f_2^{*}(x/2)+g_2^{*}(x/2)\right)^p w(x) dx \nonumber\\
& \le & C\left(\int_{{\RR}_+^2} \left(f_2^{*}(x/2)\right)^p w(x)
dx+\int_{{\RR}_+^2}
 \left(g_2^{*}(x/2)\right)^p w(x) dx\right) \nonumber\\
& \le & C\left(\int_{{\RR}_+^2} \left(f_2^{*}(x)\right)^p w(2x)
dx+\int_{{\RR}_+^2}
 \left(g_2^{*}(x)\right)^p w(2x) dx\right) \nonumber\\
& \le & C\left(\int_{{\RR}_+^2} \left(f_2^{*}(x)\right)^p w(x)
dx+\int_{{\RR}_+^2}
 \left(g_2^{*}(x)\right)^p w(x) dx\right)\nonumber\\
& = & C(\| f \|^p_{{\Lambda}^p_2(w)}+\|g\|^p_{{\Lambda}^p_2(w)}),\nonumber
\end{eqnarray}
and it follows that $\| f+g\|_{{\Lambda}^p_2(w)}\le C(\| f\|_{{\Lambda}^p_2(w)}+\|
g\|_{{\Lambda}^p_2(w)}).$

\medskip

Conversely, let $D$ and $D_1$ be two  sets of ${{\RR}^2}$ with  $D\cap
D_1=\emptyset$  and $D^{*}=D_1^{*}$, and such that if $D^{*}$  has the
representation
$$
D^{*}=\{(x_1,x_2): 0<x_1<r, 0<x_2<\phi(x_1);r>0\},
$$
(with $\phi\downarrow$), then
$$
(D\cup D_1)^{*}=\{(x_1,x_2): 0<x_1<2r, 0<x_2<\phi(x_1/2);r>0\},
$$
(this is easily done by taking $D_1$ to be a translation of the form $D_1=D+(N,0)$, where $N>0$ is
big enough). If $\|\cdot \|_{{\Lambda}^p_2(w)}$ is a quasinorm, then
$$
\| f+g \|^p_{{\Lambda}^p_2(w)}\le C( \| f \|^p_{{\Lambda}^p_2(w)}+\| g
\|^p_{{\Lambda}^p_2(w)}),
$$
and if we take $f=\chi_D$ and $g=\chi_{D_1}$, then we get

\begin{equation}\label{eqde}
\int_{(D\cup {D_1})^{*}} w(x)\, dx\leq C \int_{D^{*}}w(x)\, dx.
\end{equation}

We denote by $E:=(D\cup D_1)^{*}$, 
and by
$$
E_1:= \{(x_1,x_2): 0<x_1<2r, \phi(2 x_1)< x_2<2 \phi(x_1/2);\, r>0\}.
$$
Obviously, $ E_1\cup E= 2D^{*}$. Since $E_1^{*}=E= E^{*}$ we can apply (\ref{eqde}) with $D=E$, 
$D_1=E_1$ and  get

\begin{eqnarray}
\int_{2 D^{*}}w(x)\, dx & = & \int_{(E\cup E_1)^{*}}w(x)\, dx \le
C\int_{E^{*}}w(x)\, dx \nonumber\\
& = & C\int_{(D\cup D_1)^{*}}w(x)\, dx \le C \int_{D^{*}}w(x)\, dx\nonumber,
\end{eqnarray}
which is obviously equivalent to condition (\ref{quasi}). Thus, in view of (\ref{hom}) and (\ref{can}),
the first statement is proved.

\medskip

To prove that ${{\Lambda}^p_2(w)}$ is complete we have to show that if
$\left(f_k\right)_k\subset{{\Lambda}^p_2(w)}$
is a Cauchy sequence, then there exists a function $f\in
{{\Lambda}^p_2(w)}$ such that
 $\|f_j- f \|_{{\Lambda}^p_2(w)}\longrightarrow 0$ as $j\rightarrow\infty$.
Since $\|\cdot \|^p_{{\Lambda}^p_2(w)}$ is quasinorm and
$\left(f_k\right)_k$ is Cauchy,
 there exists a constant $C>0$ such that $\|f_j \|^p_{{\Lambda}^p_2(w)}\le
C<\infty,\,\,\forall j\in \NN$.

Also since $\left(f_j- f_k\right)_2^{*}$ is decreasing in each
variable, for a fixed $x\in\RR^2_+$, if we set $Q_x=\{y\in\RR^2_+:0<y_k\le x_k, \ k=1,2\}$,
then 
$$
 {\left(f_j- f_k\right)_2^{*}}^p (x)\int_{Q_x}w(y) dy\le
\int_{{\RR}^2_+}
 {\left(f_j- f_k\right)_2^{*}}^p (y)w(y) dy.
$$
Therefore
$$
{\left(f_j- f_k\right)_2^{*}} \longrightarrow 0, \,\,\mbox{a.e.}
$$
This implies 
$${\lambda}_{\left(f_j- f_k\right)_2^{*}}\longrightarrow 0,
\,\,\mbox{a.e.}
$$ 
and hence
$$
{\lambda}_{\left(f_j- f_k\right)}\longrightarrow 0, \,\,\mbox{a.e.},
$$
i.e., $\left(f_k\right)_k$ is Cauchy in measure. 
Hence there is a subsequence $\left(f_{k_j}\right)$ which converges
pointwise, say to a function $f$ which is measurable. 
By Proposition~\ref{propre} e) and by Fatou's lemma we have that  $f\in
{{\Lambda}^p_2(w)}$. 
Moreover,
$$
\lim_{j\rightarrow\infty} |f_{k_j}(x)-f_i(x)|=|f
(x)-f_i(x)|,\,\,x\in {{\RR}^2}.
$$
Using Fatou's lemma again and the fact that  $\left(f_k\right)_k$
is a Cauchy sequence, we finally get
$$
\|f- f_i \|_{{\Lambda}^p_2(w)}\leq C \left(\|f- f_{k_j} \|_{{\Lambda}^p_2(w)}+
\|f_i- f_{k_j} \|_{{\Lambda}^p_2(w)}\right)\longrightarrow 0,\mbox{ as }
i,j\rightarrow\infty.\quad\Box
$$

\medskip

Finally, we are now going to prove the main result of this section, namely, the characterization of
the weights $w$ for which $\Vert\cdot\Vert_{{\Lambda}^p_2(w)}$ is a norm. We begin by showing
the following necessary condition on the index $p$:

\medskip

\begin{thr}\label{indexp}
Let $0<p<\infty$.
If ${\Lambda}^p_2(w)$ is a Banach space, then $p\geq 1$.
\end{thr}

{\bf Proof.} Since ${\Lambda}^p_2(w)$ is a Banach space, there exists
$\|\cdot\|$, a norm on ${\Lambda}^p_2(w)$, such that
$$
\|f\|_{{\Lambda}^p_2(w)}\approx \|f\|.
$$
Hence
$$
\bigg\|\sum_{k=1}^{N}f_k\bigg\|_{{\Lambda}^p_2(w)}\leq C \sum_{k=1}^{N}\|f_k\| \leq
{ C} \sum_{k=1}^{N}\|f_k\|_{{\Lambda}^p_2(w)},
$$
for all $N\in {\NN}$. 
Suppose $0<p<1$ and take a decreasing sequence of domains
$$
A_{k+1}\subset A_k\subset\ldots\subset {{\RR}^2},
$$
such that $\int_{A_k^{*}} w(x)\,dx=2^{-kp}$. 
If $ f_k=2^k\chi_{A_k}$, then
$ \|f_k\|_{{\Lambda}^p_2(w)}=1$.

But for a fixed $N$,  we have that
$$
\frac{1}{N}\bigg\|\sum_{k=1}^{N}f_k\bigg\|_{{\Lambda}^p_2(w)}\leq {\tilde C}<\infty.
$$
On the other hand, since
$\left(\sum_{k=1}^{N}2^k\chi_{A_k}\right)^{*}_2=\sum_{k=1}^{N}2^k\chi_{A_k^{*}}$ (by
(\ref{reofsim})),  and
$A_{k+1}^{*}\subset A_k^{*}\subset\ldots\subset {{\RR}^2_+}$
we have (taking $A_{N+1}=\emptyset$)
\begin{eqnarray}
\frac{1}{N}\bigg\|\sum_{k=1}^{N}f_k\bigg\|_{{\Lambda}^p_2(w)} & = &
\frac{1}{N}\bigg\|\sum_{k=1}^{N}2^k\chi_{A_k}\bigg\|_{{\Lambda}^p_2(w)}
\nonumber\\
& = &  \frac{1}{N}\left(\int_{{\RR}^2_+}
\left(\sum_{k=1}^{N}2^k\chi_{A_k^{*}}\right)^p(x)w(x)dx\right)^{1/p}
 \nonumber\\
& = & \frac{1}{N}\left(\int_{{\RR}^2_+}
\left(\sum_{k=1}^{N}\left(\sum_{j=1}^{k}2^j\right)
\chi_{A_k^{*}\setminus A_{k+1}^{*}}\right)^p(x)w(x)dx\right)^{1/p}
\nonumber\\
& = & \frac{1}{N}\left(\int_{{\RR}^2_+}
\sum_{k=1}^{N}\left(\sum_{j=1}^{k}2^j\right)^p
\chi_{A_k^{*}\setminus A_{k+1}^{*}}(x)w(x)dx\right)^{1/p} \nonumber\\
& = & \frac{1}{N}\left(\sum_{k=1}^{N}\left(\sum_{j=1}^{k}2^j\right)^p
\left(\int_{A_k^{*}}
w(x)dx-\int_{A_{k+1}^{*}}(x)w(x)dx\right)\right)^{1/p} \nonumber\\
& \ge & \frac{C}{N}\left(\sum_{k=1}^{N} \left(1-
2^{-k}\right)^p\right)^{1/p}\nonumber\\
& \ge &
\frac{C}{N}\left(\sum_{k=1}^{N}2^{-p}\right)^{1/p}=C\frac{N^{1/p}}{N}\rightarrow
\infty,\,\mbox{as }N\rightarrow\infty,\nonumber
\end{eqnarray}
which is a contradiction. Hence $p\ge 1$.  $\qquad\Box$

\medskip

\begin{thr}\label{normeq}
Let $1\le p<\infty$ and $w$ be a weight in $\RR^2_+$. Then, the following conditions are equivalent:

\medskip
\noindent
a) $\Vert\cdot\Vert_{\Lambda^p_2(w)}$ is a norm.

\medskip
\noindent
b) For every  $A,B\subset\RR^2$,  $w((A\cap B)^*)+w((A\cup B)^*)\le
w(A^*)+w(B^*).$

\medskip
\noindent
c) There exists a decreasing weight $v$ on $\RR^+$ such that $w(s,t)=v(t)$,  $s,t>0$.
\end{thr}

{\bf Proof.}  
If $\Vert\cdot\Vert_{\Lambda^p_2(w)}$ is a norm, take $A,B\subset\RR^2$, $\delta>0$ and define
the functions
$$
f(x)=\left\{ \begin{array}{ll}
1+\delta, & \mbox{if $x\in A$}\\
1,  & \mbox{if $x\in (A\cup B)\setminus A$}\\
0, &  \mbox{otherwise,}
\end{array}
\right.
$$
and
$$
g(x)=\left\{ \begin{array}{ll}
1+\delta, & \mbox{if $x\in B$}\\
1,  & \mbox{if $x\in (A\cup B)\setminus B$}\\
0, &  \mbox{otherwise.}
\end{array}
\right.
$$
Then,
\begin{eqnarray*}
f^*_2(x)&=&(1+\delta)\chi_{A^*}(x)+\chi_{(A\cup B)^*\setminus A^*}(x),\\
g^*_2(x)&=&(1+\delta)\chi_{B^*}(x)+\chi_{(A\cup B)^*\setminus B^*}(x),\\
(f+g)^*_2(x)&=&(2+2\delta)\chi_{(A\cap B)^*}(x)+(2+\delta)\chi_{(A\cup B)^*\setminus
(A\cap B)^*}(x),
\end{eqnarray*}
and, hence, the triangle inequality and the fact that $1/p\le 1$ imply
\begin{eqnarray*}
\Vert f+g\Vert_{\Lambda_2^p(w)}&=&\Big((2+2\delta)^pw((A\cap B)^*)\\
&&\qquad+(2+\delta)^pw((A\cup
B)^*\setminus(A\cap B)^*)\Big)^{1/p}\\
&\le&\Vert f\Vert_{\Lambda_2^p(w)}+\Vert g\Vert_{\Lambda_2^p(w)}\\
&=&\Big((1+\delta)^pw(A^*)+w((A\cup
B)^*\setminus A^*)\Big)^{1/p}\\
&&\qquad+
\Big((1+\delta)^pw(B^*)+w((A\cup
B)^*\setminus B^*)\Big)^{1/p}\\
&\le&2^{1-1/p}\Big((1+\delta)^pw(A^*)+w((A\cup
B)^*\setminus A^*)\\
&&\qquad+(1+\delta)^pw(B^*)+w((A\cup
B)^*\setminus B^*)\Big)^{1/p}.
\end{eqnarray*}
Collecting terms, dividing both sides  by $2^{p-1}((1+\delta)^p-1)$ and letting $\delta\to0$, we
finally obtain
$$ 
w((A\cap B)^*)+w((A\cup B)^*)\le w(A^*)+w(B^*),
$$ 
which is b). Thus a) implies b).

\medskip

Assume now that b) holds. 
Fix $s,t>0$, and consider, for $\varepsilon>0$ small, the sets
\begin{eqnarray*}
A&=&(0,\varepsilon)\times(0,t)\cup (\varepsilon,s)\times(0,t-\varepsilon),\\
B&=&(0,\varepsilon)\times(0,t-\varepsilon)\cup (\varepsilon,s)\times(0,t). 
\end{eqnarray*}
Then, 
\begin{eqnarray*}
A^*&=&A,\\
B^*&=&(0,s-\varepsilon)\times(0,t)\cup (s-\varepsilon,s)\times(0,t-\varepsilon),\\
(A\cap
B)^*&=&(0,s)\times(0,t-\varepsilon),\\
(A\cup B)^*&=&(0,s)\times(0,t). 
\end{eqnarray*}
Hence using b)  we
obtain that
\begin{eqnarray*}
w((s-\varepsilon,s)\times(t-\varepsilon,t))&=&w((A\cup B)^*)-w(B^*)\le w(A^*)-w((A\cap
B)^*)\\
&=&w((0,\varepsilon)\times(t-\varepsilon,t)).
\end{eqnarray*}
Thus, dividing both sides by $\varepsilon^2$ and letting $\varepsilon\to0$ we obtain that
$w(s,t)\le w(0,t)$.

\medskip

Similarly, taking now 
\begin{eqnarray*}
A&=&(0,s)\times(0,t),\\
B&=&(0,\varepsilon)\times(\varepsilon,t+\varepsilon)\cup
(\varepsilon,s-\varepsilon)\times(0,t)\cup(s-\varepsilon,s)\times(0,t-\varepsilon), 
\end{eqnarray*}
we obtain
that 
\begin{eqnarray*}
A^*&=&A,\\
B^*&=&(0,s-\varepsilon)\times(0,t)\cup (s-\varepsilon,s)\times(0,t-\varepsilon),\\
(A\cap
B)^*&=&(0,s-2\varepsilon)\times(0,t)\cup (s-2\varepsilon,s)\times(0,t-\varepsilon),\\
(A\cup B)^*&=&(0,\varepsilon)\times(0,t+\varepsilon)\cup (\varepsilon,s)\times(0,t). 
\end{eqnarray*}
Therefore by
using b) we obtain that
\begin{eqnarray*}
w((0,\varepsilon)\times(t,t+\varepsilon))&=&w((A\cup B)^*)-w(A^*)\le w(B^*)-w((A\cap
B)^*)\\
&=&w((s-2\varepsilon,s-\varepsilon)\times(t-\varepsilon,t)).
\end{eqnarray*}
Hence, dividing both sides by $\varepsilon^2$ and letting $\varepsilon\to0$, we obtain that
$w(0,t)\le w(s,t)$ and, thus, 
$$ 
w(s,t)=w(0,t)=v(t).
$$

\medskip

To finish we will prove that $v(b)=w(0,b)\le w(0,a)=v(a)$ if $0<a\le b$: for $\varepsilon>0$ small,
take now
\begin{eqnarray*}
A&=&(0,\varepsilon)\times(0,a),\\
B&=&(0,\varepsilon)\times(\varepsilon,b).
\end{eqnarray*}
Then, 
\begin{eqnarray*}
A^*&=&A,\\
B^*&=&(0,\varepsilon)\times(0,b-\varepsilon),\\
(A\cap
B)^*&=&(0,\varepsilon)\times(0,a-\varepsilon),\\
(A\cup B)^*&=&(0,\varepsilon)\times(0,b).
\end{eqnarray*}
Hence using b) we obtain that
\begin{eqnarray*}
w((0,\varepsilon)\times(b-\varepsilon,b))&=&w((A\cup B)^*)-w(B^*)\le w(A^*)-w((A\cap
B)^*)\\
&=&w((0,\varepsilon)\times(a-\varepsilon,a)).
\end{eqnarray*}
 Thus, dividing both sides by $\varepsilon^2$ and letting $\varepsilon\to0$ we obtain that
$w(0,b)\le w(0,a)$.

\medskip

Finally, we are now going to prove that  c) implies a).
By Theorem~\ref{iterrea}  we know  that $f^*_2(s,t)=(f_x^{*y}(t))^{*x}(s)$. Thus, using
the fact that $\Vert\cdot\Vert_{\Lambda^p(v)}$ is a norm, if $v$ is decreasing (see \cite{L}), and
Minkowski's inequality, we obtain 

\begin{eqnarray*}
\Vert f+g\Vert_{\luw}&=&\bigg(\int_{\Rdm}[(f+g)^*_2(s,t)]^pw(s,t)\,dsdt\bigg)^{1/p}\\
&=&\bigg(\int_0^{\infty}\bigg(\int_0^{\infty}\Big[\Big((f_x+g_x)^{*y}(t)\Big)^{*x}(s)
\Big]^p\,ds\bigg)v(t)\,dt\bigg)^{1/p}\\
&=&\bigg(\int_0^{\infty}\bigg(\int_{\RR}\Big[(f_x+g_x)^{*y}(t)\Big]^p\,dx\bigg)v(t)\,dt
\bigg)^{1/p}\\ 
&=&\bigg(\int_{\RR}
\bigg(\int_0^{\infty}\Big[(f_x+g_x)^{*y}(t)\Big]^pv(t)\,dt\bigg)\,dx\bigg)^{1/p}\\
&\le&\bigg(\int_{\RR}\bigg[
\bigg(\int_0^{\infty}\Big[(f_x)^{*y}(t)\Big]^pv(t)\,dt\bigg)^{1/p}\\
&&\qquad+
\bigg(\int_0^{\infty}\Big[(g_x)^{*y}(t)\Big]^pv(t)\,dt\bigg)^{1/p}\bigg]^p\,dx\bigg)^{1/p}\\
&\le&\bigg(\int_{\RR}
\bigg(\int_0^{\infty}\Big[(f_x)^{*y}(t)\Big]^pv(t)\,dt\bigg)\,dx\bigg)^{1/p}\\
&&\qquad+
\bigg(\int_{\RR}
\bigg(\int_0^{\infty}\Big[(g_x)^{*y}(t)\Big]^pv(t)\,dt\bigg)\,dx\bigg)^{1/p}\\
&=&\bigg(\int_{\Rdm}(f^*_2(s,t))^pw(s,t)\,dsdt\bigg)^{1/p}\\
&&\qquad+
\bigg(\int_{\Rdm}(g^*_2(s,t))^pw(s,t)\,dsdt\bigg)^{1/p}\\
&=&\Vert
f\Vert_{\luw}+\Vert g\Vert_{\luw}. 
\end{eqnarray*}
Thus, in view of (\ref{hom}) and (\ref{can}), the proof is complete. $\qquad\Box$
\medskip

\begin{rem}\rm
Observe that the equivalences proved in Theorem~\ref{normeq} in particular say that
$\luw=L^p(\Lambda^p(v,dy),dx)$, which is a mixed norm space.
\end{rem}

\vskip 1cm

\bigbreak


\begin{thebibliography}{BPSo}

\bibitem[AM]{AM} M.A. Ari\~no and  B. Muckenhoupt, {\sl Maximal functions on classical Lorentz spaces and
Hardy's  inequality with weights for nonincreasing functions}, Trans. Amer. Math. Soc. {\bf 320} (1990),
727--735.

\bibitem[Ba]{Ba} S.\ Barza, {\sl Weighted multidimensional integral inequalities and applications},
Ph.D. Thesis, Lule\aa\ University, 1999.

\bibitem[BPSo]{BPSo} S.\ Barza, L.E.\ Persson, and J.\ Soria, {\sl Sharp weighted
multidimensional integral inequalities for monotone functions}, Math. Nachr. {\bf 210} (2000),
43--58.

\bibitem[BPSt]{BPSt} S.\ Barza, L.\,E.\ Persson, and V.\ Stepanov, {\sl On weighted
multidimensional embeddings for monotone functions}, Math. Scand. {\bf 88} (2001),
303--319.

\bibitem[BS]{BS} C. Bennet and R. Sharpley, {\sl Interpolation of Operators}, Academic Press, 1988.


\bibitem[Bl]{Bl} A. P. \ Blozinski, {\sl Multivariate rearrangements and Banach function spaces with
mixed norms}, Trans. Amer. Math. Soc. {\bf 263} (1981), 149--167.

\bibitem[CS]{CS} M.J. Carro and J. Soria,  {\sl Weighted Lorentz spaces and the Hardy operator}, J. Funct. Anal. {\bf 112}
(1993), 480--494.

\bibitem[LL]{LL} E.H.\ Lieb and M.\  Loss, {\sl Analysis}, Graduate Studies in
Mathematics, 14. American Mathematical Society, Providence, RI, 2001.

\bibitem[L]{L} G.G. Lorentz,  {\sl On the theory of spaces $\Lambda$}, Pacific J.
Math. {\bf 1} (1951), 411--429.

\bibitem[Sa]{Sa} E. Sawyer,  {\sl Boundedness of classical operators on classical
Lorentz spaces}, Studia Math. {\bf 96} (1990), 145--158.

\bibitem[St]{St} V. Stepanov,  {\sl The weighted Hardy's inequality for nonincreasing functions},
Trans. Amer. Math. Soc. {\bf 338} (1993), 173--186.


\end{thebibliography}
\end{document}